\documentclass[12pt,reqno]{amsart}
\usepackage[T1]{fontenc}
\usepackage{lmodern,microtype}
\usepackage[margin=1in]{geometry}
\usepackage{amsmath,amssymb,amsthm,mathtools}
\usepackage[colorlinks=true,linkcolor=black,citecolor=black,urlcolor=blue]{hyperref}
\allowdisplaybreaks[2]
\newtheorem{theorem}{Theorem}[section]
\newtheorem{lemma}[theorem]{Lemma}
\newtheorem{proposition}[theorem]{Proposition}

\theoremstyle{definition}

\newtheorem{example}[theorem]{Example}
\theoremstyle{remark}
\newtheorem{remark}[theorem]{Remark}
\numberwithin{equation}{section}

\newcommand{\C}{\mathbb C}
\newcommand{\R}{\mathbb R}
\newcommand{\Irr}{\operatorname{Irr}}
\newcommand{\ad}{\operatorname{ad}}
\newcommand{\id}{\operatorname{id}}
\newcommand{\distq}{\operatorname{dist}_{q}}

\begin{document}
\title[Quantum GH convergence for discrete quantum groups]{Quantum Gromov--Hausdorff Convergence of Spectral Truncations for Amenable Discrete Quantum Groups with Modified Rapid Decay}
\author{Xintao Peng}
\begin{abstract}
Let $G$ be an amenable discrete quantum group with the modified rapid decay
property relative to a proper central length $\ell$. We construct a
modularly corrected commutator seminorm on $C_r^*(G)$ and prove its right
invariance. Averaging over the coproduct yields a regular bi-invariant
Lip-norm. For the spectral projections $P_R=\chi_{[0,R]}(D_\ell)$, we
equip the finite-dimensional operator systems
$A_R=P_RC_r^*(G)P_R$ with the Lip-norms induced by the compressed left and
right coactions. We prove quantum Gromov--Hausdorff convergence and obtain
an explicit approximation estimate. These spectral compressions are
length-filtered Peter--Weyl truncations of the compact dual. The resulting
regular bi-invariant Lip-norm also places the family within Leimbach's
operator Gromov--Hausdorff convergence theorem. The construction applies
to non-unimodular discrete quantum groups, including discrete duals of
standard deformations of compact semisimple Lie groups.
\end{abstract}
\maketitle

\section{Introduction}

Compact quantum metric spaces provide a noncommutative
counterpart of compact metric spaces~\cite{RieffelCQMS}.
Their finite-dimensional approximations have been studied
via Fourier, Peter--Weyl and spectral truncations; see,
for example, \cite{ConnesSuijlekom,RieffelFourier,
LeimbachSuijlekom,GaudillotEstrada,Leimbach}.
In particular, Leimbach~\cite{Leimbach} established
convergence of Peter--Weyl compressions for coamenable
compact quantum groups under a regular bi-invariant
Lip-norm hypothesis. In our previous work
~\cite{PengWangSpectral}, we studied spectral truncations
of compact quantum groups and proved their convergence
in quantum Gromov--Hausdorff distance.

We now turn to amenable discrete quantum groups satisfying
the modified rapid decay property. Our aim is to construct
suitable Lip-norms and study the convergence of the
associated spectral truncations. The modular correction
in the RD estimate plays an essential role, particularly
for non-unimodular quantum groups. We use the construction
of Bhowmick, Voigt and Zacharias~\cite{BVZ}. Given a proper
central length $\ell$, let $D_\ell$ be the associated
length operator. Set
\[
P_R=\chi_{[0,R]}(D_\ell),
\qquad
A_R=P_R C_r^*(G)P_R.
\]
By properness of $\ell$, each $A_R$ is a finite-dimensional
operator system.

It is useful to make explicit the relation with Peter--Weyl
truncations. If $\widehat G$ denotes the compact dual of $G$ and
$p_\alpha$ is the projection onto the coefficient space associated
with $\alpha\in\operatorname{Irr}(\widehat G)$, then
\[
 L^2(G)=\bigoplus_{\alpha\in\operatorname{Irr}(\widehat G)}H_\alpha,
 \qquad
 D_\ell=\sum_{\alpha}\ell(\alpha)p_\alpha.
\]
Consequently,
\[
 P_R=\sum_{\ell(\alpha)\leq R}p_\alpha.
\]
Thus $A_R$ is precisely the compression associated with the
length-filtered finite set
\[
 \Lambda_R=\{\alpha\in\operatorname{Irr}(\widehat G):
                  \ell(\alpha)\leq R\}.
\]
In particular, these are length-filtered Peter--Weyl truncations; we
do not regard spectral truncation as a different class of compression.
The point here is to construct the metric structure required for this
particular family directly from the modified RD estimates on $G$.

More precisely, the modular correction is described by a positive
self-adjoint operator $T$, which is generally unbounded and whose
eigenvalues depend on the row index of an irreducible corepresentation.
The associated seminorm
\[
 L_T^k(a)=
 \big\|T^k\operatorname{ad}(D_\ell)^k(a)T^k\big\|,
 \qquad a\in\mathbb C[G]_{\rm sa},
\]
is right-invariant. Since the modular weights need not be scalar on
irreducible blocks, we do not assume that this original seminorm is
bi-invariant. Instead, we form its bi-invariant envelope
\[
 \widetilde L_T^k(a)=
 \sup_{\omega\in S(A)}
 L_T^k\big((\omega\otimes\operatorname{id})
                  \widehat\Delta(a)\big).
\]
We then use both compressed coactions to define
\[
 L_R(x)=\max\{L_R^r(x),L_R^\ell(x)\}.
\]
The two terms are needed separately: the compression map is controlled
by the bi-invariance of $\widetilde L_T^k$, whereas the symbol map is
controlled by the left coaction.

Our contribution is the construction and verification of this
modularly corrected regular bi-invariant Lip-norm from modified RD,
together with the resulting quantitative approximation estimates for
the length-filtered compressions. This is not a new abstract
Peter--Weyl convergence mechanism. Rather, it supplies a concrete
class of Lip-norms on compact duals of amenable discrete quantum groups,
including non-unimodular examples, to which the general convergence
theory of Leimbach applies. We also give a direct order-unit
quantum Gromov--Hausdorff estimate.

The paper is organized as follows. Section~2 recalls the modified RD
estimates and the order-unit approximation criterion. Section~3
constructs the invariant seminorms and the compressed coactions, and
proves convergence. Section~4 records examples.

\section{Preliminaries}
\subsection{Discrete quantum groups, lengths and RD}
Let $G$ be a discrete quantum group with compact dual $\widehat G$. Write
\[
 C_0(G)=\bigoplus_{\alpha\in\Irr(\widehat G)}^{c_0}B(H_\alpha),
 \qquad A=C_r^*(G)=C(\widehat G),
\]
with coproducts $\Delta$ and $\widehat\Delta$, respectively. The canonical dense Hopf $*$-algebra of $A$ is denoted by $\C[G]$; the algebraic direct sum of the matrix blocks of $C_0(G)$ is denoted by $C_c(G)$. The Fourier transform gives inverse linear isomorphisms
\[
 \mathcal F:C_c(G)\longrightarrow\C[G],\qquad
 \mathcal G:\C[G]\longrightarrow C_c(G).
\]
The left Haar weight on $C_c(G)$ gives the Hilbert space $L^2(G)$. We use
its Fourier identification with the GNS space of the Haar state of $A$ and
follow the conventions of~\cite{BVZ,Vergnioux}. The quotient map
$C_f^*(G)\to A$ is an isomorphism precisely when $G$ is amenable. In this
case the algebraic counit extends to a state $\epsilon$ on $A$.

A length function is a map $\ell:\Irr(\widehat G)\to[0,\infty)$ such that $\ell(\mathbf1)=0$, $\ell(\bar\alpha)=\ell(\alpha)$, and $\ell(\gamma)\leq\ell(\alpha)+\ell(\beta)$ whenever $\gamma$ is contained in $\alpha\otimes\beta$. It is proper if each set $\{\alpha:\ell(\alpha)\leq R\}$ is finite and $\ell(\alpha)=0$ only for the trivial corepresentation. The corresponding central length is $L=\sum_\alpha\ell(\alpha)p_\alpha$. Properness is precisely the condition ensuring that each bounded spectral subspace of the associated length operator is finite dimensional. In particular, the spectral projections $P_R$ used below have finite rank and increase strongly to the identity.

Let $\widehat h$ be the Haar state of $A$ and let $\widehat\Lambda:A\to L^2(G)$ be its GNS map. The associated length operator is
\[
 D_\ell=\bigoplus_{\alpha\in\Irr(\widehat G)}\ell(\alpha)I_{H_\alpha},
\]
where $H_\alpha$ is the coefficient subspace spanned by $\widehat\Lambda(u_{ij}^\alpha)$. Properness makes its spectral projections finite rank.

Let $F=\bigoplus_\alpha F^\alpha$ be the modular element of $G$.
Each $F^\alpha$ is positive and invertible, with
\[
 \operatorname{Tr}(F^\alpha)
 =\operatorname{Tr}((F^\alpha)^{-1})
 =\dim_q(\alpha).
\]
Choose an orthonormal basis which diagonalizes every $F^\alpha$, and
write
\[
 F^\alpha e_i^\alpha=f_i^\alpha e_i^\alpha,
 \qquad f_i^\alpha>0.
\]
Put
\[
 \xi_{ij}^\alpha=\widehat\Lambda(u_{ij}^\alpha),
 \qquad
 H_{\alpha,i}
 =\overline{\operatorname{span}}
   \{\xi_{ij}^\alpha:1\leq j\leq d_\alpha\}.
\]
The GNS Hilbert space decomposes orthogonally as
\[
 L^2(G)=\bigoplus_{\alpha\in\Irr(\widehat G)}
             \bigoplus_{i=1}^{d_\alpha}H_{\alpha,i}.
\]
Set
\[
 c_{\alpha,i}
 =\frac{\dim_q(\alpha)}{d_\alpha}f_i^\alpha,
 \qquad d_\alpha=\dim(\alpha).
\]
On the coefficient space associated with $\alpha$, the modular
correction is the matrix
\[
 C^\alpha=\frac{\dim_q(\alpha)}{d_\alpha}F^\alpha.
\]
Define the positive self-adjoint direct-sum operators $F,C,T$ by
\[
 F|_{H_{\alpha,i}}=f_i^\alpha I,\qquad
 C|_{H_{\alpha,i}}=c_{\alpha,i}I,\qquad
 T=C^{1/(2k)}|_{H_{\alpha,i}}
   =c_{\alpha,i}^{1/(2k)}I.
\]
Their domains are determined by the corresponding square-summability
conditions; explicitly,
\[
 \begin{aligned}
 \operatorname{Dom}(F)
 &=\left\{\xi=\bigoplus_{\alpha,i}\xi_{\alpha,i}:
       \sum_{\alpha,i}(f_i^\alpha)^2
       \|\xi_{\alpha,i}\|^2<\infty\right\},\\
 \operatorname{Dom}(C)
 &=\left\{\xi=\bigoplus_{\alpha,i}\xi_{\alpha,i}:
       \sum_{\alpha,i}c_{\alpha,i}^{\,2}
       \|\xi_{\alpha,i}\|^2<\infty\right\},\\
 \operatorname{Dom}(T)
 &=\left\{\xi=\bigoplus_{\alpha,i}\xi_{\alpha,i}:
       \sum_{\alpha,i}c_{\alpha,i}^{\,1/k}
       \|\xi_{\alpha,i}\|^2<\infty\right\}.
 \end{aligned}
\]
These operators are positive and injective, and may be unbounded.
Their finite-support vectors form a common invariant core, and
$T^k=C^{1/2}$. Since $D_\ell$ acts as
$\ell(\alpha)I$ on each $H_{\alpha,i}$, the operators $F,C,T$ and
$D_\ell$ strongly commute. The unimodular case is
characterized by $F^\alpha=I$ and hence $C=T=I$.

For $x=\sum_{\alpha,i,j}x_{ij}^\alpha u_{ij}^\alpha\in\C[G]$, the
modified Sobolev norm of~\cite{BVZ} is
\begin{equation}\label{eq:sobolev}
 \|x\|_{2,s}^2
 =\sum_{\alpha,i,j}
   \frac{(1+\ell(\alpha))^{2s}}{d_\alpha}
   |x_{ij}^\alpha|^2.
\end{equation}
Equivalently, it is the norm of
$(1+L)^s\mathcal G(C^{1/2}*x)$ in the left-Haar GNS space.
Thus, in the non-unimodular case, it is not the uncorrected
Haar-GNS norm of $x$. We say that $G$ has the modified property RD
with respect to $\ell$ if there are constants $c,s>0$ such that
\begin{equation}\label{eq:RD}
 \|\mathcal F(f)\|_A\leq c\|f\|_{2,s},
 \qquad f\in C_c(G).
\end{equation}
By the Fourier isometry, this is equivalent to
\[
 \|x\|_A\leq c\|x\|_{2,s},
 \qquad x\in\C[G].
\]

For $a\in\C[G]$, define
\[
 \ad(D_\ell)^0(a)=a,\qquad
 \ad(D_\ell)^{j+1}(a)
 =[D_\ell,\ad(D_\ell)^j(a)].
\]
The modularly corrected commutator seminorm is
\begin{equation}\label{eq:LT}
 L_T^k(a)=
 \big\|T^k\ad(D_\ell)^k(a)T^k\big\|.
\end{equation}
The product in~\eqref{eq:LT} is initially defined on the common
finite-support core. By the estimates in~\cite[Lemma~6.5]{BVZ},
it extends to a bounded operator for $a\in\C[G]$. We use the same
notation for this bounded extension. No boundedness of $C$ or $T$
is assumed.

\begin{theorem}[Proper-length form of the BVZ criterion]\label{thm:BVZ}
Suppose $G$ has a proper central length $\ell$ and satisfies~\eqref{eq:RD}. If $k>s$ is a positive integer, then $L_T^k$, with domain $\C[G]_{\rm sa}$, is a Lip-norm on $A_{\rm sa}$ in the sense that its Monge--Kantorovich metric induces the weak$^*$ topology on $S(A)$.
\end{theorem}
\begin{proof}
By~\cite[Lemma~6.5]{BVZ}, $L_T^k$ is a Lipschitz seminorm with scalar kernel. The coefficient estimate of~\cite[Lemma~6.6]{BVZ} gives
\begin{equation}\label{eq:coeff}
 \sum_{\alpha,i,j}\frac{\ell(\alpha)^{2k}}{d_\alpha}|a_{ij}^\alpha|^2
 \leq L_T^k(a)^2.
\end{equation}
For $a_{>N}=\sum_{\ell(\alpha)>N,i,j}a_{ij}^\alpha u_{ij}^\alpha$, the RD estimate and~\eqref{eq:sobolev} imply, for $N\geq1$,
\begin{align}
 \|a_{>N}\|&\leq c\|a_{>N}\|_{2,s}\notag\\
 &\leq c\,2^sN^{s-k}
 \left(\sum_{\alpha,i,j}\frac{\ell(\alpha)^{2k}}{d_\alpha}|a_{ij}^\alpha|^2\right)^{1/2}\notag\\
 &\leq c\,2^sN^{s-k}L_T^k(a).\label{eq:tail}
\end{align}
Fix the Haar state $\widehat h$ and consider the centered unit ball
\[
 E=\{a\in\C[G]_{\rm sa}:L_T^k(a)\leq1,\ \widehat h(a)=0\}.
\]
This is the centering required in the Ozawa--Rieffel criterion. Write
\[
 a_{\leq N}=\sum_{\ell(\alpha)\leq N,i,j}a_{ij}^\alpha u_{ij}^\alpha,
 \qquad a_{>N}=a-a_{\leq N}.
\]
For fixed $N$, properness implies that only finitely many irreducibles occur in $a_{\leq N}$. Since $\widehat h(a)=0$, the coefficient of the trivial corepresentation vanishes. For every remaining coefficient with $\ell(\alpha)\leq N$, the number $\ell(\alpha)^{2k}/d_\alpha$ is strictly positive and fixed. Thus~\eqref{eq:coeff} bounds every coefficient of $a_{\leq N}$ uniformly for $a\in E$. It follows that
\[
 \{a_{\leq N}:a\in E\}
\]
is bounded in a finite-dimensional vector space and hence is totally bounded in the $C^*$-norm. On the other hand,~\eqref{eq:tail} gives
\[
 \sup_{a\in E}\|a-a_{\leq N}\|\leq c\,2^sN^{s-k}\longrightarrow0.
\]
A finite net for the low-frequency truncations, enlarged by this uniform tail bound, is therefore a finite net for $E$. Hence $E$ is totally bounded. The Ozawa--Rieffel criterion~\cite[Proposition~6.2]{BVZ} proves the Lip-norm property.

Bhowmick, Voigt and Zacharias state~\cite[Theorem~7.4]{BVZ}
for finitely generated discrete quantum groups equipped with word
length. The argument above uses only the finiteness of bounded
length balls, which is exactly the properness assumption. Thus the
same low-frequency and high-frequency estimates establish the
Lip-norm property for an arbitrary proper central length.
\end{proof}

\subsection{Compact quantum metric spaces}
Let $X$ be a unital $C^*$-algebra and $X_{\rm sa}$ its self-adjoint part.
A seminorm $L$ on a dense real subspace of $X_{\rm sa}$ is a Lip-norm if
$L(1_X)=0$, $\ker L=\mathbb R1_X$, and the metric
\[
 d_L(\mu,\nu)=\sup\{ |\mu(x)-\nu(x)|:x\in X_{\rm sa},\ L(x)\leq1\}
\]
induces the weak$^*$ topology on $S(X)$. We write
\[
 \operatorname{diam}(X,L)=\sup_{\mu,\nu\in S(X)}d_L(\mu,\nu).
\]
We use Rieffel's equivalent criterion: for one (hence every) state $\omega$,
\[
 \{x\in X_{\rm sa}:L(x)\leq1,\ \omega(x)=0\}
\]
is norm totally bounded~\cite{RieffelCQMS}. The same definitions apply to the self-adjoint order-unit space of an operator system. In particular,
$A_R$ is treated as an operator system, not as a $C^*$-subalgebra; its
seminorms are defined on $(A_R)_{\rm sa}$ and have kernel $\mathbb RP_R$.

For compact quantum metric order-unit spaces, an admissible Lip-norm on
$A\oplus B$ is a Lip-norm whose coordinate quotient seminorms are $L_A$
and $L_B$. The quantum Gromov--Hausdorff distance is the infimum of the
Hausdorff distances between the two state spaces in the metric induced by
such an admissible Lip-norm; see also~\cite{KaadKyed} for background on
 quantum metric structures. We recall an order-unit approximation criterion for quantum Gromov--Hausdorff convergence, a convenient formulation of the standard contractive approximation argument; see, e.g., \cite{KaadKyed,Suijlekom}.

\begin{proposition}[Contractive approximation criterion]\label{prop:GHcriterion}
Let $(A,L_A)$ and $(B,L_B)$ be compact quantum metric order-unit spaces.
Suppose unital positive maps $\Phi:A\to B$ and $\Psi:B\to A$ satisfy
\[
 L_B(\Phi(a))\leq L_A(a),\qquad
 L_A(\Psi(b))\leq L_B(b),
\]
and, for some $\eta\geq0$,
\[
 \|\Psi\Phi(a)-a\|\leq\eta L_A(a),\qquad
 \|\Phi\Psi(b)-b\|\leq\eta L_B(b)
\]
on the self-adjoint Lip-domains. Then
\[
 \distq((A,L_A),(B,L_B))\leq\eta.
\]
\end{proposition}
\begin{proof}
Fix $r>\eta$ and put
\[
 M_r(a,b)=\max\{L_A(a),L_B(b),r^{-1}\|b-\Phi(a)\|,
                         r^{-1}\|a-\Psi(b)\|\}.
\]
The quotient seminorms are $L_A$ and $L_B$: the lower bounds follow from
the first two terms, while the reverse bounds follow by testing
$(a,\Phi(a))$ and $(\Psi(b),b)$ and using $\eta<r$. The kernel of $M_r$
is the common scalar line.

Choose states $\omega_A\in S(A)$ and $\omega_B\in S(B)$. If $M_r(a,b)\leq1$,
then the centered elements
\[
 a^\circ=a-\omega_A(a)1_A,\qquad
 b^\circ=b-\omega_B(b)1_B
\]
lie in the respective centered Lip-unit balls, which are totally bounded.
Moreover,
\[
\begin{aligned}
 |\omega_B(b)-\omega_A(a)|
 &\leq \|b-\Phi(a)\|
       +|\omega_B(\Phi(a))-\omega_A(a)|\\
 &\leq r+d_{L_A}(\omega_B\circ\Phi,\omega_A)\\
 &\leq r+\operatorname{diam}(A,L_A).
\end{aligned}
\]
Thus, after subtracting the common scalar $\omega_A(a)$, the pairs in the
$M_r$-unit ball lie in a totally bounded subset of
$A_{\rm sa}\oplus B_{\rm sa}$. Hence $M_r$ is an admissible Lip-norm.
For $\mu\in S(A)$, pair $\mu$ with $\mu\circ\Psi\in S(B)$. Whenever
$M_r(a,b)\leq1$,
\[
 |\mu(a)-(\mu\circ\Psi)(b)|\leq\|a-\Psi(b)\|\leq r.
\]
The analogous estimate pairs each state of $B$ with its pullback through
$\Phi$. The Hausdorff distance is at most $r$; letting $r\downarrow\eta$
proves the claim.
\end{proof}

\subsection{Spectral compressions and quantum distance}
Let $\lambda:A\to B(L^2(G))$ be the left regular representation, and
suppress $\lambda$ when no confusion can arise. For $R>0$, set
\[
\begin{aligned}
 P_R&=\chi_{[0,R]}(D_\ell),&
 \mathcal H_R&=P_RL^2(G),\\
 A_R&=P_R\lambda(A)P_R\subseteq B(\mathcal H_R),&
 \tau_R(a)&=P_R\lambda(a)P_R.
\end{aligned}
\]
\begin{lemma}\label{lem:compression}
For every $R>0$, the projection $P_R$ has finite rank,
\[
 \operatorname{rank}(P_R)=\sum_{\ell(\alpha)\leq R}d_\alpha^2,
\]
and $A_R$ is a finite-dimensional operator system with order unit $P_R$. The compression $\tau_R:A\to A_R$ is ucp and surjective. Moreover,
\[
 \tau_R(\C[G])=A_R.
\]
\end{lemma}
\begin{proof}
Only finitely many irreducibles satisfy $\ell(\alpha)\leq R$. Since the
coefficient space of $u^\alpha$ has dimension $d_\alpha^2$,
\[
 \operatorname{rank}(P_R)=\sum_{\ell(\alpha)\leq R}d_\alpha^2.
\]
Thus $A_R$ is finite dimensional. Compression is completely positive and
$\tau_R(1_A)=P_R$, so $\tau_R:A\to A_R$ is ucp and is onto by the
definition of $A_R$. Density of $\C[G]$ in $A$ implies that
$\tau_R(\C[G])$ is dense in $A_R$; being a linear subspace of a
finite-dimensional space, it is closed. Hence $\tau_R(\C[G])=A_R$.
\end{proof}

We use the order-unit version of Rieffel's quantum Gromov--Hausdorff distance described above. The contractive approximation criterion in Proposition~\ref{prop:GHcriterion} is the only form needed below.

\section{Spectral truncation convergence}
Throughout this section $G$ is amenable, and hence the counit $\epsilon$
is a state on $A$. Let $W,V$ be the standard implementing unitaries on
$L^2(G)\otimes L^2(G)$, so that
\[
 \widehat\Delta(a)=W^*(1\otimes a)W=V(a\otimes1)V^*.
\]

\subsection{Invariance of the modularly corrected seminorm}
Choose bases so that $F^\alpha=\operatorname{diag}(f_1^\alpha,\ldots,f_{d_\alpha}^\alpha)$. Then
\[
 T\widehat\Lambda(u_{ij}^\alpha)=\tau_i^\alpha\widehat\Lambda(u_{ij}^\alpha),\qquad
 \tau_i^\alpha=\left(\frac{\dim_q(\alpha)}{d_\alpha}f_i^\alpha\right)^{1/(2k)}.
\]
The dependence on the row index is relevant when $G$ is non-unimodular.

\begin{lemma}\label{lem:unitary-commutation}
On the algebraic core $\widehat\Lambda(\C[G])\odot\widehat\Lambda(\C[G])$, the unitaries $W,W^*,V,V^*$ preserve the core and strongly commute with the spectral projections of $1\otimes D_\ell$ (for $W,W^*$) and $D_\ell\otimes1$ (for $V,V^*$). Moreover,
\[
 [1\otimes D_\ell,W]=[1\otimes D_\ell,W^*]=0,
 \quad [D_\ell\otimes1,V]=[D_\ell\otimes1,V^*]=0,
 \quad [T\otimes1,V]=[T\otimes1,V^*]=0.
\]
The commutation relations with $T\otimes1$ are strong commutation relations.
\end{lemma}
\begin{proof}
Put $\xi_{ij}^\alpha=\widehat\Lambda(u_{ij}^\alpha)$ and
\[
 H_{\alpha,i}=\overline{\operatorname{span}}
 \{\xi_{ij}^\alpha:1\leq j\leq d_\alpha\},
 \qquad
 H_\beta=\overline{\operatorname{span}}
 \{\xi_{pq}^\beta:1\leq p,q\leq d_\beta\}.
\]
With our convention for the implementing unitaries, the algebraic formulas
are
\begin{align*}
 W^*(\xi_{ij}^\alpha\otimes\xi_{pq}^\beta)
 &=\sum_r\widehat\Lambda(u_{pr}^\beta u_{ij}^\alpha)\otimes\xi_{rq}^\beta,\\
 V(\xi_{ij}^\alpha\otimes\xi_{pq}^\beta)
 &=\sum_r\xi_{ir}^\alpha\otimes\widehat\Lambda(u_{rj}^\alpha u_{pq}^\beta).
\end{align*}
Products of matrix coefficients are finite linear combinations of
coefficients. The corresponding formulas for the inverse unitaries give
\[
 W^{\pm1}(L^2(G)\otimes H_\beta)=L^2(G)\otimes H_\beta,
 \qquad
 V^{\pm1}(H_{\alpha,i}\otimes L^2(G))=H_{\alpha,i}\otimes L^2(G).
\]
In particular, the relevant subspaces are reducing. Since
$D_\ell|_{H_\beta}=\ell(\beta)I$ and
$D_\ell|_{H_{\alpha,i}}=\ell(\alpha)I$, we obtain
\[
 [1\otimes D_\ell,W]=[1\otimes D_\ell,W^*]=0,
 \qquad
 [D_\ell\otimes1,V]=[D_\ell\otimes1,V^*]=0.
\]
Moreover, $T|_{H_{\alpha,i}}=\tau_i^\alpha I$, so the same reducing
subspaces give
\[
 [T\otimes1,V]=[T\otimes1,V^*]=0.
\]
These are strong commutation relations: the unitaries commute with the
spectral projections of the indicated self-adjoint operators.
\end{proof}

\begin{remark}
For non-Kac $G$, one cannot in general infer $[1\otimes T,W]=0$:
the second row index changes from $p$ to $r$, and the corresponding
$T$-eigenvalues may differ. We therefore use only the right invariance of
$L_T^k$ and obtain the other invariance by averaging in~\eqref{eq:averaged-L}.
\end{remark}

\begin{proposition}\label{prop:right-invariance}
For $a\in\C[G]$ and $\varphi\in S(A)$,
\[
 L_T^k((\id\otimes\varphi)\widehat\Delta(a))\leq L_T^k(a).
\]
Thus $L_T^k$ is right-invariant.
\end{proposition}
\begin{proof}
Lemma~\ref{lem:unitary-commutation} gives, on the algebraic core,
\[
 \ad(D_\ell\otimes1)^k(\widehat\Delta(a))
 =V(\ad(D_\ell)^k(a)\otimes1)V^*.
\]
The strong commutation with $T\otimes1$ gives
\[
 (T^k\otimes1)\ad(D_\ell\otimes1)^k(\widehat\Delta(a))(T^k\otimes1)
 =V(T^k\ad(D_\ell)^k(a)T^k\otimes1)V^*.
\]
Write $\widehat\Delta(a)=\sum_{j=1}^m a_j\otimes b_j$ with
$a_j,b_j\in\C[G]$. The weighted commutator is then
\[
 \sum_{j=1}^m T^k\ad(D_\ell)^k(a_j)T^k\otimes b_j
 \in B(L^2(G))\odot A,
\]
so the slice $\id\otimes\varphi$ is defined and contractive. Since it
acts on the second leg, it commutes with the first-leg commutators and with
multiplication by $T^k$. Consequently,
\[
\begin{aligned}
 L_T^k((\id\otimes\varphi)\widehat\Delta(a))
 &=\|(\id\otimes\varphi)(V(T^k\ad(D_\ell)^k(a)T^k\otimes1)V^*)\|\\
 &\leq\|\varphi\|\,\|T^k\ad(D_\ell)^k(a)T^k\|\\
 &=L_T^k(a),
\end{aligned}
\]
where the slice map is contractive because $\varphi$ is a state.
\end{proof}

Define, for $a\in\C[G]_{\rm sa}$,
\begin{equation}\label{eq:averaged-L}
 \widetilde L_T^k(a)=\sup_{\omega\in S(A)}L_T^k((\omega\otimes\id)\widehat\Delta(a)).
\end{equation}

\begin{theorem}\label{thm:bi-invariant}
The seminorm $\widetilde L_T^k$ is finite on $\C[G]_{\rm sa}$, is bi-invariant, and is a Lip-norm on $A$.
\end{theorem}

\begin{proof}
For $a\in\C[G]$, write the finite algebraic coproduct as
\[
\widehat\Delta(a)=\sum_{j=1}^m a_j\otimes b_j
\]
with $a_j,b_j\in\C[G]$. For every state $\omega$,
\[
(\omega\otimes\id)\widehat\Delta(a)
=\sum_j\omega(a_j)b_j\in\C[G].
\]
The triangle inequality for $L_T^k$ gives
\[
\widetilde L_T^k(a)
\leq \sum_j\|a_j\|L_T^k(b_j)<\infty.
\]
Thus $\widetilde L_T^k$ is a seminorm. Since the counit is a state and
$(\epsilon\otimes\id)\widehat\Delta=\id$, we have
\begin{equation}\label{eq:comparison-L}
L_T^k(a)\leq\widetilde L_T^k(a).
\end{equation}

For $\omega,\mu\in S(A)$, put
\[
\alpha_\omega=(\omega\otimes\id)\widehat\Delta,
\qquad
\beta_\mu=(\id\otimes\mu)\widehat\Delta.
\]
Coassociativity gives
\[
\alpha_\omega\alpha_\mu=\alpha_{\mu*\omega},
\qquad
\alpha_\omega\beta_\mu=\beta_\mu\alpha_\omega,
\]
where
\[
\mu*\omega=(\mu\otimes\omega)\widehat\Delta.
\]
Consequently,
\[
\begin{aligned}
\widetilde L_T^k(\alpha_\mu(a))
&=\sup_{\omega\in S(A)}
L_T^k(\alpha_{\mu*\omega}(a))\\
&\leq\sup_{\nu\in S(A)}
L_T^k(\alpha_\nu(a))\\
&=\widetilde L_T^k(a),
\end{aligned}
\]
since $\mu*\omega\in S(A)$ for every $\omega\in S(A)$.
On the other hand, right invariance of $L_T^k$ gives
\[
\begin{aligned}
\widetilde L_T^k(\beta_\mu(a))
&=\sup_{\omega\in S(A)}
L_T^k(\alpha_\omega\beta_\mu(a))\\
&=\sup_{\omega\in S(A)}
L_T^k(\beta_\mu\alpha_\omega(a))\\
&\leq\sup_{\omega\in S(A)}
L_T^k(\alpha_\omega(a))\\
&=\widetilde L_T^k(a).
\end{aligned}
\]
Thus $\widetilde L_T^k$ is bi-invariant.

Moreover,~\eqref{eq:comparison-L} gives
\[
\ker\widetilde L_T^k
\subseteq\ker L_T^k=\mathbb R1.
\]
Since $\widetilde L_T^k(1)=0$, its kernel is exactly $\mathbb R1$.

Fix a state $\omega_0$ and consider the centered unit ball
\[
E_{\widetilde L}
=
\left\{
a\in\C[G]_{\rm sa}:
\widetilde L_T^k(a)\leq1,\ 
\omega_0(a)=0
\right\}.
\]
By~\eqref{eq:comparison-L}, this set is contained in the corresponding centered unit ball for $L_T^k$, which is norm totally bounded by Theorem~\ref{thm:BVZ}. Rieffel's criterion therefore proves that $\widetilde L_T^k$ is a Lip-norm.
\end{proof}

\subsection{The coactions on the compressed operator systems}
We use the standard operator-system definition: a right coaction is a unital completely isometric map $\alpha:X\to X\otimes A$ satisfying coassociativity and
$\overline{\operatorname{span}}\{(1_X\otimes A)\alpha(X)\}=X\otimes A$; the left-handed definition is analogous. Two coactions cocommute when the two iterated maps agree.

All tensor products below are the concrete minimal operator-system tensor
products. The fixed-point space of a right coaction is
\[
 X^\alpha=\{x\in X:\alpha(x)=x\otimes1\};
\]
for a left coaction $\beta:X\to A\otimes X$ it is defined by $\beta(x)=1\otimes x$. Ergodicity means that the relevant fixed-point space consists only of scalar multiples of the order unit.

\begin{proposition}[Compressed coactions]
\label{thm:coactions}
For every $R>0$, the formulas
\[
 \alpha_R(x)=V(x\otimes1)V^*,\qquad
 \beta_R(x)=W^*(1\otimes x)W,\qquad x\in A_R,
\]
define completely isometric right and left coactions
\[
 \alpha_R:A_R\to A_R\otimes A,\qquad
 \beta_R:A_R\to A\otimes A_R.
\]
They satisfy the Podle\'s density conditions, are ergodic, and cocommute:
\begin{equation}\label{eq:cocommute}
 (\beta_R\otimes\id)\alpha_R=(\id\otimes\alpha_R)\beta_R.
\end{equation}
\end{proposition}
\begin{proof}
The projection $P_R$ is the Peter--Weyl projection associated with
$\Lambda_R=\{\alpha:\ell(\alpha)\leq R\}$. The general compression
construction gives the unital completely isometric coactions
$\alpha_R$ and $\beta_R$, their coassociativity, and the Podle\'s
density conditions; see~\cite[Proposition~7.4]{Leimbach}. We verify
the identities and additional properties needed here.

By Lemma~\ref{lem:unitary-commutation} and the spectral theorem,
\[
 [P_R\otimes1,V]=0,\qquad [1\otimes P_R,W]=0.
\]
Hence the formulas for $\alpha_R$ and $\beta_R$ take values in
$A_R\otimes A$ and $A\otimes A_R$, respectively. For $a\in A$,
\begin{equation}\label{eq:coaction-compatibility}
 \alpha_R(\tau_R(a))
 =(\tau_R\otimes\id)\widehat\Delta(a),\qquad
 \beta_R(\tau_R(a))
 =(\id\otimes\tau_R)\widehat\Delta(a).
\end{equation}
These identities are independent of the choice of lift because
$\tau_R$ is surjective and the compressed coproducts are well-defined.

Let $\widehat h$ be the Haar state of $A$. If
$\alpha_R(x)=x\otimes1$, write $x=\tau_R(a)$. Applying
$\id\otimes\widehat h$ to~\eqref{eq:coaction-compatibility} gives
\[
 x=(\id\otimes\widehat h)\alpha_R(x)
 =\tau_R((\id\otimes\widehat h)\widehat\Delta(a))
 =\widehat h(a)P_R.
\]
Thus $A_R^{\alpha_R}=\C P_R$. The same argument, applying
$\widehat h\otimes\id$, gives $A_R^{\beta_R}=\C P_R$.

Finally, for $x=\tau_R(a)$, compatibility gives
\[
 \begin{aligned}
 (\beta_R\otimes\id)\alpha_R(x)
 &=(\id\otimes\tau_R\otimes\id)
   (\widehat\Delta\otimes\id)\widehat\Delta(a),\\
 (\id\otimes\alpha_R)\beta_R(x)
 &=(\id\otimes\tau_R\otimes\id)
   (\id\otimes\widehat\Delta)\widehat\Delta(a).
 \end{aligned}
\]
The two expressions agree by coassociativity of
$\widehat\Delta$. This proves ergodicity and cocommutativity.
\end{proof}

\begin{remark}
The existence of the compressed coactions is part of the general
Peter--Weyl compression framework of~\cite{Leimbach}. The points used
specifically below are the compatibility identities
\eqref{eq:coaction-compatibility}, the ergodicity of the two
coactions, and their cocommutativity for the length-filtered family
$(P_R)_{R>0}$. Properness ensures that $P_R$ has finite rank and
that $\tau_R(\C[G])=A_R$.
\end{remark}

\begin{lemma}\label{lem:counit}
Let $\epsilon$ be the bounded counit of $A$. If $\gamma:X\to X\otimes A$ is an injective right coaction, then
\[
 (\id_X\otimes\epsilon)\gamma=\id_X.
\]
If $\eta:X\to A\otimes X$ is an injective left coaction, then
\[
 (\epsilon\otimes\id_X)\eta=\id_X.
\]
\end{lemma}
\begin{proof}
For a right coaction, coassociativity gives
\[
 \gamma((\id_X\otimes\epsilon)\gamma(x))
 =(\id_X\otimes\id_A\otimes\epsilon)
   (\id_X\otimes\widehat\Delta)\gamma(x)
 =\gamma(x),
\]
where the last equality uses $(\id_A\otimes\epsilon)\widehat\Delta=\id_A$. Injectivity of $\gamma$ proves the first identity. For a left coaction, the same argument, using
$(\epsilon\otimes\id_A)\widehat\Delta=\id_A$, proves the second.
\end{proof}

\subsection{The truncation Lip-norm and convergence}
For $x\in(A_R)_{\rm sa}$ define
\begin{align}
 L_R^r(x)&=\sup_{\psi\in S(A_R)}\widetilde L_T^k((\psi\otimes\id)\alpha_R(x)),\notag\\
 L_R^\ell(x)&=\sup_{\phi\in S(A_R)}\widetilde L_T^k((\id\otimes\phi)\beta_R(x)),\label{eq:LR}\\
 L_R(x)&=\max\{L_R^r(x),L_R^\ell(x)\}.\notag
\end{align}
For $\phi\in S(A_R)$, let
\begin{equation}\label{eq:symbol}
 \rho_R^\phi(x)=(\id\otimes\phi)\beta_R(x),\qquad x\in A_R.
\end{equation}

\begin{proposition}\label{prop:LR}
$L_R$ is a Lip-norm on the finite-dimensional operator system $A_R$. Moreover, for $a\in\C[G]_{\rm sa}$ and $x\in(A_R)_{\rm sa}$,
\begin{equation}\label{eq:contractions}
 L_R(\tau_R(a))\leq\widetilde L_T^k(a),\qquad
 \widetilde L_T^k(\rho_R^\phi(x))\leq L_R(x).
\end{equation}
\end{proposition}
\begin{proof}
Choose $a\in\C[G]$ with $x=\tau_R(a)$ and write
$\widehat\Delta(a)=\sum_{j=1}^m a_j\otimes b_j$ with
$a_j,b_j\in\C[G]$. By~\eqref{eq:coaction-compatibility},
\[
 (\psi\otimes\id)\alpha_R(x)
 =\sum_{j=1}^m\psi(\tau_R(a_j))b_j,
 \qquad
 (\id\otimes\phi)\beta_R(x)
 =\sum_{j=1}^m\phi(\tau_R(b_j))a_j.
\]
Consequently,
\[
\begin{aligned}
 L_R^r(x)&\leq\sum_{j=1}^m\|a_j\|\,\widetilde L_T^k(b_j),\\
 L_R^\ell(x)&\leq\sum_{j=1}^m\|b_j\|\,\widetilde L_T^k(a_j).
\end{aligned}
\]
In particular, both suprema in~\eqref{eq:LR} are finite. If $L_R^r(x)=0$, then
\[
 \widetilde L_T^k((\psi\otimes\id)\alpha_R(x))=0
 \quad\text{for all }\psi\in S(A_R).
\]
Since $\ker\widetilde L_T^k=\R1$, every such slice is scalar. States separate points of $A_R$, so $\alpha_R(x)\in A_R\otimes\C1$. The right counit identity from Lemma~\ref{lem:counit} gives $\alpha_R(x)=x\otimes1$, and ergodicity gives $x\in\R P_R$. Conversely, unitality shows that every scalar multiple of $P_R$ lies in the kernel. The same argument for $\beta_R$ gives
\[
 \ker L_R^r=\ker L_R^\ell=\R P_R,
 \qquad \ker L_R=\R P_R.
\]
A finite-dimensional seminorm with scalar kernel is a Lip-norm, so the first assertion follows.

For $a\in\C[G]_{\rm sa}$, the compatibility identities for the coactions and the bi-invariance of $\widetilde L_T^k$ give
\[
 L_R^r(\tau_R(a))\leq\widetilde L_T^k(a),\qquad
 L_R^\ell(\tau_R(a))\leq\widetilde L_T^k(a).
\]
Taking the maximum proves the first inequality in~\eqref{eq:contractions}. The second follows immediately from the definition of $L_R^\ell$ and~\eqref{eq:symbol}. The symbol map is ucp, as it is a composition of a coaction and a state slice.
\end{proof}

\begin{lemma}\label{lem:slice}
For $a\in\C[G]_{\rm sa}$, $x\in(A_R)_{\rm sa}$ and $\mu,\nu\in S(A)$,
\begin{align}
 \| (\id\otimes\mu)\widehat\Delta(a)-(\id\otimes\nu)\widehat\Delta(a)\|
 &\leq d_{\widetilde L_T^k}(\mu,\nu)\widetilde L_T^k(a),\label{eq:sliceA}\\
 \| (\id\otimes\mu)\alpha_R(x)-(\id\otimes\nu)\alpha_R(x)\|
 &\leq d_{\widetilde L_T^k}(\mu,\nu)L_R^r(x).\label{eq:sliceR}
\end{align}
\end{lemma}
\begin{proof}
For~\eqref{eq:sliceA}, set
$z=(\id\otimes\mu)\widehat\Delta(a)-(\id\otimes\nu)\widehat\Delta(a)$.
Since $z=z^*$,
\[
 \|z\|=\sup_{\omega\in S(A)}|\omega(z)|.
\]
For each $\omega$, put $y_\omega=(\omega\otimes\id)\widehat\Delta(a)$. Then $y_\omega=y_\omega^*$ and left invariance gives $\widetilde L_T^k(y_\omega)\leq\widetilde L_T^k(a)$. Hence
\[
 |\omega(z)|=|\mu(y_\omega)-\nu(y_\omega)|
 \leq d_{\widetilde L_T^k}(\mu,\nu)\widetilde L_T^k(a).
\]
Taking the supremum gives~\eqref{eq:sliceA}. For~\eqref{eq:sliceR},
put $z_R=(\id\otimes\mu)\alpha_R(x)-(\id\otimes\nu)\alpha_R(x)$.
For each $\psi\in S(A_R)$, let
$y_\psi=(\psi\otimes\id)\alpha_R(x)$. Then
$\widetilde L_T^k(y_\psi)\leq L_R^r(x)$, and the same argument, now
using $\|z_R\|=\sup_{\psi\in S(A_R)}|\psi(z_R)|$, proves
~\eqref{eq:sliceR}.
\end{proof}

\begin{theorem}[Quantum Gromov--Hausdorff convergence]\label{thm:main}
Let $G$ be an amenable discrete quantum group with modified property RD for a proper central length $\ell$, and choose an integer $k>s$ as in Theorem~\ref{thm:BVZ}. Equip $A=C_r^*(G)$ with $\widetilde L_T^k$ and $A_R=P_RAP_R$ with $L_R$. Then
\[
 \distq\big((A,\widetilde L_T^k),(A_R,L_R)\big)\longrightarrow0
 \qquad(R\to\infty).
\]
\end{theorem}
\begin{proof}
We first construct states $\phi_R\in S(A_R)$ such that
\[
 \mu_R:=\phi_R\circ\tau_R\longrightarrow\epsilon
 \quad\text{in the weak$^*$ topology.}
\]
Since $\lambda$ is faithful, the convex hull of its vector states is
weak$^*$ dense in $S(A)$. The metric $d_{\widetilde L_T^k}$ induces the
weak$^*$ topology. Hence, for each $n$, choose a finite convex combination
$\omega^{(n)}$ of vector states with
\[
 d_{\widetilde L_T^k}(\omega^{(n)},\epsilon)<\frac1{2n}.
\]
Write the vectors implementing this convex combination as $\xi_1,\ldots,\xi_m$. Since $P_R\to I$ strongly, the normalized vectors
\[
 \xi_j^{(R)}=\frac{P_R\xi_j}{\|P_R\xi_j\|}
\]
are defined for all sufficiently large $R$ and converge to $\xi_j$. The associated vector states converge weak$^*$, hence in
$d_{\widetilde L_T^k}$, to the original vector states. Their finite convex
combination $\omega_R^{(n)}$ therefore converges to $\omega^{(n)}$ in the
same metric. More explicitly, if $\omega_R^{(n)}$ is implemented by
$\xi_1^{(R)},\ldots,\xi_m^{(R)}$ with weights $t_1,\ldots,t_m$, then
\[
 \omega_R^{(n)}(a)=\sum_{j=1}^m t_j
 \langle\lambda(a)\xi_j^{(R)},\xi_j^{(R)}\rangle
 =\phi_R^{(n)}(\tau_R(a))
\]
for a state $\phi_R^{(n)}$ on $A_R$. Choose an increasing sequence
$R_n\to\infty$ such that
\[
 d_{\widetilde L_T^k}(\omega_R^{(n)},\omega^{(n)})<\frac1{2n}
 \quad\text{whenever }R\geq R_n.
\]
For $R_n\leq R<R_{n+1}$, set $\phi_R=\phi_R^{(n)}$; the choice for
$R<R_1$ is immaterial. Then
\[
 d_R:=d_{\widetilde L_T^k}(\phi_R\circ\tau_R,\epsilon)<\frac1n
 \quad(R_n\leq R<R_{n+1}),
\]
so $d_R\to0$.
Put $\rho_R=\rho_R^{\phi_R}$. By~\eqref{eq:contractions}, $\tau_R$ and $\rho_R$ are unital positive Lip-contractions.

For $a\in\C[G]_{\rm sa}$, coaction compatibility gives
\[
 \rho_R\tau_R(a)=(\id\otimes\mu_R)\widehat\Delta(a).
\]
Since $a=(\id\otimes\epsilon)\widehat\Delta(a)$,~\eqref{eq:sliceA} yields
\[
 \|\rho_R\tau_R(a)-a\|\leq d_R\widetilde L_T^k(a).
\]
For $x\in(A_R)_{\rm sa}$, the two induced coactions satisfy
$(\tau_R\otimes\id)\beta_R=(\id\otimes\tau_R)\alpha_R$, and hence
\[
 \tau_R\rho_R(x)=(\id\otimes\mu_R)\alpha_R(x).
\]
Using the right counit identity and~\eqref{eq:sliceR}, we obtain
\[
 \|\tau_R\rho_R(x)-x\|\leq d_R L_R^r(x)\leq d_RL_R(x).
\]
Apply Proposition~\ref{prop:GHcriterion} with $\eta=d_R$. It follows that
\[
 \distq\big((A,\widetilde L_T^k),(A_R,L_R)\big)\leq d_R\to0.
\]
\end{proof}

\begin{remark}[Comparison with Leimbach's operator-Gromov--Hausdorff theorem]
The convergence mechanism in Theorem~\ref{thm:main} is closely
related to, but should be distinguished from, the general
Peter--Weyl convergence theorem of Leimbach~\cite[Theorem~7.13]{Leimbach}.
We verify that its hypotheses are satisfied for the family at hand.

First, amenability of $G$ is equivalent to coamenability of its
compact dual $\widehat G$. Properness of $\ell$ implies that
\[
 \Irr(\widehat G)
 =\bigcup_{n=1}^\infty
   \{\alpha:\ell(\alpha)\leq n\},
\]
a countable union of finite sets. Hence $A=C(\widehat G)$ is
separable. Moreover, if $R\leq S$, then
\[
 P_R\leq P_S,\qquad P_R\vee P_S=P_S,
\]
and
\[
 P_R\longrightarrow1
 \quad\text{strongly as }R\longrightarrow\infty.
\]
Thus $(P_R)_{R>0}$ is an increasing Peter--Weyl net of the type
required in~\cite[Theorem~7.13]{Leimbach}.

Second, Theorem~\ref{thm:bi-invariant} shows that
$\widetilde L_T^k$ is a bi-invariant Lip-norm. It is regular in
the sense used there, since its domain contains the dense Hopf
$*$-algebra $\C[G]$. Finally, the seminorm $L_R$ in
\eqref{eq:LR} is precisely the maximum of the seminorms induced
by the compressed right and left coactions, which is the
operator-system Lip-norm appearing in~\cite{Leimbach}.

Consequently, all hypotheses of~\cite[Theorem~7.13]{Leimbach}
are satisfied, and we obtain the stronger convergence
\[
 \operatorname{dist}_{\mathrm{op}}
 \big((A_R,L_R),(A,\widetilde L_T^k)\big)
 \longrightarrow0.
\]
The role of the present paper is to construct the required
regular bi-invariant Lip-norm from the modified RD property of
an amenable discrete quantum group, with the modular correction
retained in the non-unimodular case. We also give a
proper-length version of the BVZ Lip-norm criterion, so that
finite generation is not imposed on the discrete quantum group.
The abstract Peter--Weyl convergence principle itself is not
claimed as a new result here. The direct estimates in
Theorem~\ref{thm:main} additionally give the explicit
order-unit quantum Gromov--Hausdorff bound $d_R$.

\end{remark}

\section{Examples}
\begin{example}[Groups of polynomial growth]
Let $\Gamma$ be a finitely generated discrete group of polynomial growth, with a proper word length $\ell$. Such a group is amenable and has RD; if its growth degree is $d$, one may take any RD exponent $s>d/2$. Since $\Gamma$ is unimodular, $T=1$. Thus the theorem applies to $A=C_r^*(\Gamma)$, including $\Gamma=\mathbb Z^n$ and the discrete Heisenberg group
\[
 H_3(\mathbb Z)=\left\{\begin{pmatrix}1&a&c\\0&1&b\\0&0&1\end{pmatrix}:a,b,c\in\mathbb Z\right\},
\]
which has polynomial growth of degree four and admits an RD exponent $s>2$.
\end{example}

\begin{example}[Non-Kac examples]
Let $K$ be a simply connected compact semisimple Lie group and $K_q$ its standard deformation, $0<q<1$. By~\cite[Proposition~4.3]{BVZ}, the discrete dual $G=\widehat{K_q}$ is amenable and has the modified property RD. Since $K_q$ is non-Kac for $q\ne1$, this gives non-unimodular examples covered by the theorem.

For $K_q=SU_q(2)$ the irreducibles are indexed by $n\in\mathbb N_0$, with ordinary dimension $n+1$, quantum dimension $[n+1]_q$, and standard proper length $\ell(n)=n$. In a suitable orthonormal basis the modular matrix is diagonal with eigenvalues $q^{-n},q^{-n+2},\ldots,q^n$ (up to reversing the basis order). Thus
\[
 C_n=\frac{[n+1]_q}{n+1}\operatorname{diag}(q^{-n},q^{-n+2},\ldots,q^n),
 \qquad T_n=C_n^{1/(2k)}.
\]
The modular correction is nontrivial for $q\ne1$, so this example tests the genuinely non-Kac part of the construction. The same argument applies to $SU_q(3)$ and, more generally, to the standard deformations covered by~\cite[Proposition~4.3]{BVZ}.
\end{example}

\begin{remark}
The group examples above are unimodular, so $F=C=T=1$ and the construction reduces to the higher commutator seminorm associated with the length operator. The compact semisimple $q$-deformations show that the argument is not restricted to the tracial case: the row-dependent modular weights are retained in the right-invariance calculation and are handled by the two-sided truncation seminorm.
\end{remark}


\begin{thebibliography}{20}
\bibitem{BVZ} J. Bhowmick, C. Voigt and J. Zacharias, Compact quantum metric spaces from quantum groups of rapid decay, \textit{J. Noncommut. Geom.} \textbf{9} (2015), no.~4, 1175--1200.
\bibitem{ConnesSuijlekom} A. Connes and W. D. van Suijlekom, Spectral truncations in noncommutative geometry and operator systems, \textit{Commun. Math. Phys.} \textbf{383} (2021), 2021--2067.
\bibitem{GaudillotEstrada} Y. Gaudillot-Estrada and W. D. van Suijlekom, Convergence of spectral truncations for compact metric groups, \textit{Int. Math. Res. Not. IMRN} (2025), no.~13, rnaf197.
\bibitem{KaadKyed} J. Kaad and D. Kyed, \textit{The Quantum Metric Structure of Quantum $SU(2)$}, Memoirs of the European Mathematical Society, vol.~18, EMS Press, 2025.
\bibitem{Leimbach} M. Leimbach, Convergence of Peter--Weyl truncations of compact quantum groups, \textit{J. Noncommut. Geom.} (2025), doi:10.4171/JNCG/634.
\bibitem{LeimbachSuijlekom} M. Leimbach and W. D. van Suijlekom, Gromov--Hausdorff convergence of spectral truncations for tori, \textit{Adv. Math.} \textbf{439} (2024), 109496.
\bibitem{Li} H. Li, Compact quantum metric spaces and ergodic actions of compact quantum groups, \textit{J. Funct. Anal.} \textbf{256} (2009), no.~10, 3368--3408.
\bibitem{PengWangSpectral}
X. Peng and Q. Wang,
\newblock \emph{Gromov--Hausdorff Convergence of Spectral
Truncations for Quantum Groups},
\newblock arXiv:2604.27475 [math.OA], 2026.
\bibitem{RieffelCQMS} M. A. Rieffel, Compact quantum metric spaces, in \textit{Operator Algebras, Quantization, and Noncommutative Geometry}, Contemp. Math., vol.~365, Amer. Math. Soc., 2004, pp.~315--330.
\bibitem{RieffelFourier} M. A. Rieffel, Convergence of Fourier truncations for compact quantum groups and finitely generated groups, \textit{J. Geom. Phys.} \textbf{192} (2023), 104921.
\bibitem{Suijlekom}
W. D. van Suijlekom,
Gromov--Hausdorff convergence of state spaces for spectral truncations,
\textit{J. Geom. Phys.} \textbf{162} (2021), 104075.

\bibitem{Vergnioux} R. Vergnioux, The property of rapid decay for discrete quantum groups, \textit{J. Operator Theory} \textbf{57} (2007), no.~2, 303--324.

\end{thebibliography}
\end{document}